# Optimality

## The Second Erich L. Lehmann Symposium

Javier Rojo, Editor

# Contents

## Copulas and Decoupling

## Regression Trees

## Competing Risks

## Robustness

## Multiscale Stochastic Processes

## Asymptotics

## Density Estimation

# Brief history of the Lehmann Symposia: Origins, goals and motivation

The idea of the Lehmann Symposia as platforms to encourage a revival of interest in fundamental questions in theoretical statistics, while keeping in focus issues that arise in contemporary interdisciplinary cutting-edge scientific problems, developed during a conversation that I had with Victor Perez Abreu during one of my visits to Centro de Investigación en Matemáticas (CIMAT) in Guanajuato, Mexico. Our goal was and has been to showcase relevant theoretical work to encourage young researchers and students to engage in such work.

The First Lehmann Symposium on Optimality took place in May of 2002 at Centro de Investigación en Matemáticas in Guanajuato, Mexico. A brief account of the Symposium has appeared in Vol. 44 of the Institute of Mathematical Statistics series of Lecture Notes and Monographs. The volume also contains several works presented during the First Lehmann Symposium. All papers were refereed. The program and a picture of the participants can be found on-line at the website `http://www.stat.rice.edu/lehmann/lst-Lehmann.html`.

The Second Lehmann Symposium on Optimality was held from May 19–May 22, 2004 at Rice University. There were close to 175 participants in the Symposium. A partial list and a photograph of participants, as well as the details of the scientific program, are provided in the next few pages. All scientific activities took place in Duncan Hall in the School of Engineering. Most of the plenary and invited speakers agreed to be videotaped and their talks may be accessed by visiting the following website: `http://webcast.rice.edu/webcast.php?action=details&event=408`. All papers presented in this volume were refereed, and one third of submitted papers were rejected.

At the time of this writing, plans are underway to hold the Third Lehmann Symposium at the Mathematical Sciences Research Institute during May of 2007.

I want to acknowledge the help from members of the Scientific Program Committee: Jane-Ling Wang (UC Davis), David W. Scott (Rice University), Juliet P. Shaffer (UC Berkeley), Deborah Mayo (Virginia Polytechnic Institute), Jef Teugels (Katholieke Universiteit Leuven), James R. Thompson (Rice University), and Javier Rojo (Chair).

The Symposia could not take place without generous financial support from various institutions. The First Symposium was financed in its entirety by CIMAT under the direction of Victor Perez Abreu. The Second Lehmann Symposium was generously funded by The National Science Foundation, Pfizer, Inc., The University of Texas MD Anderson Cancer Center, CIMAT, and Cytel. Shulamith Gross at NSF, Demissie Alemayehu at Pfizer, Gary Rosner at MD Anderson Cancer Center, Victor Perez Abreu at CIMAT, and Cyrrus Mehta at Cytel, encouraged and facilitated the process to obtain the support. The Rice University School of Engineering's wonderful physical facilities were made available for the Symposium at no charge.

Finally, thanks go the Statistics Department at Rice University for facilitating my participation in these activities.

May 15th, 2006
Javier Rojo
Rice University
Editor

# Contributors to this volume

Aaberge, R., *Statistics Norway*

Baraniuk, R. G., *Rice University*
Basu, A. K., *Calcutta University*
Bhattacharya, D., *Visva-Bharati University*
Bjerve, S., *University of Oslo*

Cheng, C., *St. Jude Children's Research Hospital*
Cox, D. R., *Nuffield College, Oxford*

Dabrowska, D. M., *University of California, Los Angeles*
de la Peña, V. H., *Columbia University*
Doksum, K., *University of Wisconsin, Madison*

El Barmi, H., *Baruch College, City University of New York*

Freidlin, B., *National Cancer Institute*

Gastwirth, J. L., *The George Washington University*

Ibragimov, R., *Harvard University*
Ibrahim, J., *University of North Carolina*

Juárez, S., *Veracruzana University*

Leeb, H., *Yale University*
Lehmann, E. L., *University of California, Berkeley*
Loh, W.-Y., *University of Wisconsin, Madison*

Mayo, D. G., *Virginia Polytechnical Institute*
Mnatsakanov, R. M., *West Virginia University*
Mukerjee, H., *Wichita State University*

Ribeiro, V. J., *Rice University*
Riedi, R. H., *Rice University*
Romano, J. P., *Stanford University*
Ruymgaart, F. H., *Texas Tech University*

Schucany, W. R., *Southern Methodist University*
Shaffer, J. P., *University of California*
Shaikh, A. M., *Stanford University*
Sharakhmetov, S., *Tashkent State Economics University*
Singpurwalla, N. D., *The George Washington University*
Spanos, A., *Virginia Polytechnical Institute and State University*
Székely, G. J., *Bowling Green State University, Hungarian Academy of Sciences*

Yin, G., *MD Anderson Cancer Center*

Zheng, G., *National Heart, Lung and Blood Institute*

## SCIENTIFIC PROGRAM

### The Second Erich L. Lehmann Symposium
### May 19–22, 2004
### Rice University

**Symposium Chair and Organizer Javier Rojo**

Statistics Department, MS-138
Rice University
6100 Main Street
Houston, TX 77005

**Co-Chair Victor Perez-Abreu**

Probability and Statistics
CIMAT
Callejon Jalisco S/N
Guanajuato, Mexico

### Plenary Speakers

| | |
|---|---|
| **Erich L. Lehmann**<br>UC Berkeley | *Conflicting principles in hypothesis testing* |
| **Peter Bickel**<br>UC Berkeley | *From rank tests to semiparametrics* |
| **Ingram Olkin**<br>Stanford University | *Probability models for survival and reliability analysis* |
| **D. R. Cox**<br>Nuffield College<br>Oxford | *Graphical Markov models: A tool for interpretation* |
| **Emanuel Parzen**<br>Texas A&M University | *Data modeling, quantile/quartile functions, confidence intervals, introductory statistics reform* |
| **Bradley Efron**<br>Stanford University | *Confidence regions and inferences for a multivariate normal mean vector* |
| **Kjell Doksum**<br>UC Berkeley and<br>UW Madison | *Modeling money* |
| **Persi Diaconis**<br>Stanford University | *In praise of statistical theory* |

# Invited Sessions

## *New Investigators*

| | |
|---|---|
| **Javier Rojo,** | Organizer |
| **William C. Wojciechowski,** | Chair |
| **Gabriel Huerta**<br>U of New Mexico | *Spatio-temporal analysis of Mexico city ozone levels* |
| **Sergio Juarez**<br>U Veracruzana Mexico | *Robust and efficient estimation for the generalized Pareto distribution* |
| **William C. Wojciechowski**<br>Rice University | *Adaptive robust estimation by simulation* |
| **Rudolf H. Riedi**<br>Rice University | *Optimal sampling strategies for tree-based time series* |

## *Multiple hypothesis tests: New approaches—optimality issues*

| | |
|---|---|
| **Juliet P. Shaffer,** | Chair |
| **Juliet P. Shaffer**<br>UC Berkeley | *Different types of optimality in multiple testing* |
| **Joseph Romano**<br>Stanford University | *Optimality in stepwise hypothesis testing* |
| **Peter Westfall**<br>Texas Tech University | *Optimality considerations in testing massive numbers of hypotheses* |

## *Robustness*

| | |
|---|---|
| **James R. Thompson,** | Chair |
| **Adrian Raftery**<br>U of Washington | *Probabilistic weather forecasting using Bayesian model averaging* |
| **James R. Thompson**<br>Rice University | *The simugram: A robust measure of market risk* |
| **Nozer D. Singpurwalla**<br>George Washington U | *The hazard potential: An approach for specifying models of survival* |

## *Extremes and Finance*

| | |
|---|---|
| **Jef Teugels**, | Chair |
| **Richard A. Davis**<br>Colorado State University | *Regular variation and financial time series models* |
| **Hansjoerg Albrecher**<br>University of Graz<br>Austria | *Ruin theory in the presence of dependent claims* |
| **Patrick L. Brockett**<br>U of Texas, Austin | *A chance constrained programming approach to pension plan management when asset returns are heavy tailed* |

### *Recent Advances in Longitudinal Data Analysis*

**Naisyin Wang**, Chair

**Raymond J. Carroll**, Texas A&M Univ. — *Semiparametric efficiency in longitudinal marginal models*

**Pushing Hsieh**, UC Davis — *Some issues and results on nonparametric maximum likelihood estimation in a joint model for survival and longitudinal data*

**Jane-Ling Wang**, UC Davis — *Functional regression and principal components analysis for sparse longitudinal data*

### *Semiparametric and Nonparametric Testing*

**David W. Scott,** Chair

**Jeffrey D. Hart**, Texas A&M Univ. — *Semiparametric Bayesian and frequentist tests of trend for a large collection of variable stars*

**Joseph Gastwirth**, George Washington U. — *Efficiency robust tests for linkage or association*

**Irene Gijbels**, U Catholique de Louvain — *Nonparametric testing for monotonicity of a hazard rate*

### *Philosophy of Statistics*

**Persi Diaconis,** Chair

**David Freedman**, UC Berkeley — *Some reflections on the foundations of statistics*

**Sir David Cox**, Nuffield College, Oxford — *Some remarks on statistical inference*

**Deborah Mayo**, Virginia Tech — *The theory of statistics as the "frequentist's" theory of inductive inference*

## Special contributed session

**Shulamith T. Gross,** Chair

**Victor Hugo de la Pena**, Columbia University — *Pseudo maximization and self-normalized processes*

**Wei-Yin Loh**, U of Wisconsin, Madison — *Regression tree models for data from designed experiments*

**Shulamith T. Gross**, NSF and Baruch College/CUNY — *Optimizing your chances of being funded by the NSF*

## Contributed papers

**Aris Spanos,** Virginia Tech: *Where do statistical models come from? Revisiting the problem of specification*

**Hannes Leeb,** Yale University: *The large-sample minimal coverage probability of confidence intervals in regression after model selection*

**Jun Yan,** University of Iowa: *Parametric inference of recurrent alternating event data*

**Gâbor J. Székely,** Bowling Green State U and Hungarian Academy of Sciences: *Student's t-test for scale mixture errors*

**Jaechoul Lee,** Boise State University: *Periodic time series models for United States extreme temperature trends*

**Loki Natarajan,** University of California, San Diego: *Estimation of spontaneous mutation rates*

**Chris Ding,** Lawrence Berkeley Laboratory: *Scaled principal components and correspondence analysis: clustering and ordering*

**Mark D. Rothmann,** Biologies Therapeutic Statistical Staff, CDER, FDA: *Inferences about a life distribution by sampling from the ages and from the obituaries*

**Victor de Oliveira,** University of Arkansas: *Bayesian inference and prediction of Gaussian random fields based on censored data*

**Jose Aimer T. Sanqui,** Appalachian State University: *The skew-normal approximation to the binomial distribution*

**Guosheng Yin,** The University of Texas MD Anderson Cancer Center: *A class of Bayesian shared gamma frailty models with multivariate failure time data*

**Eun-Joo Lee,** Texas Tech University: *An application of the Hâjek–Le Cam convolution theorem*

**Daren B. H. Cline,** Texas A&M University: *Determining the parameter space, Lyapounov exponents and existence of moments for threshold ARCH and GARCH time series*

**Hammou El Barmi,** Baruch College: *Restricted estimation of the cumulative incidence functions corresponding to $K$ competing risks*

**Asheber Abebe,** Auburn University: *Generalized signed-rank estimation for nonlinear models*

**Yichuan Zhao,** Georgia State University: *Inference for mean residual life and proportional mean residual life model via empirical likelihood*

**Cheng Cheng,** St. Jude Children's Research Hospital: *A significance threshold criterion for large-scale multiple tests*

**Yuan-Ji,** The University of Texas MD Anderson Cancer Center: *Bayesian mixture models for complex high-dimensional count data*

**K. Krishnamoorthy,** University of Louisiana at Lafayette: *Inferences based on generalized variable approach*

**Vladislav Karguine,** Cornerstone Research: *On the Chernoff bound for efficiency of quantum hypothesis testing*

**Robert Mnatsakanov,** West Virginia University: *Asymptotic properties of moment-density and moment-type CDF estimators in the models with weighted observations*

**Bernard Omolo,** Texas Tech University: *An aligned rank test for a repeated observations model with orthonormal design*

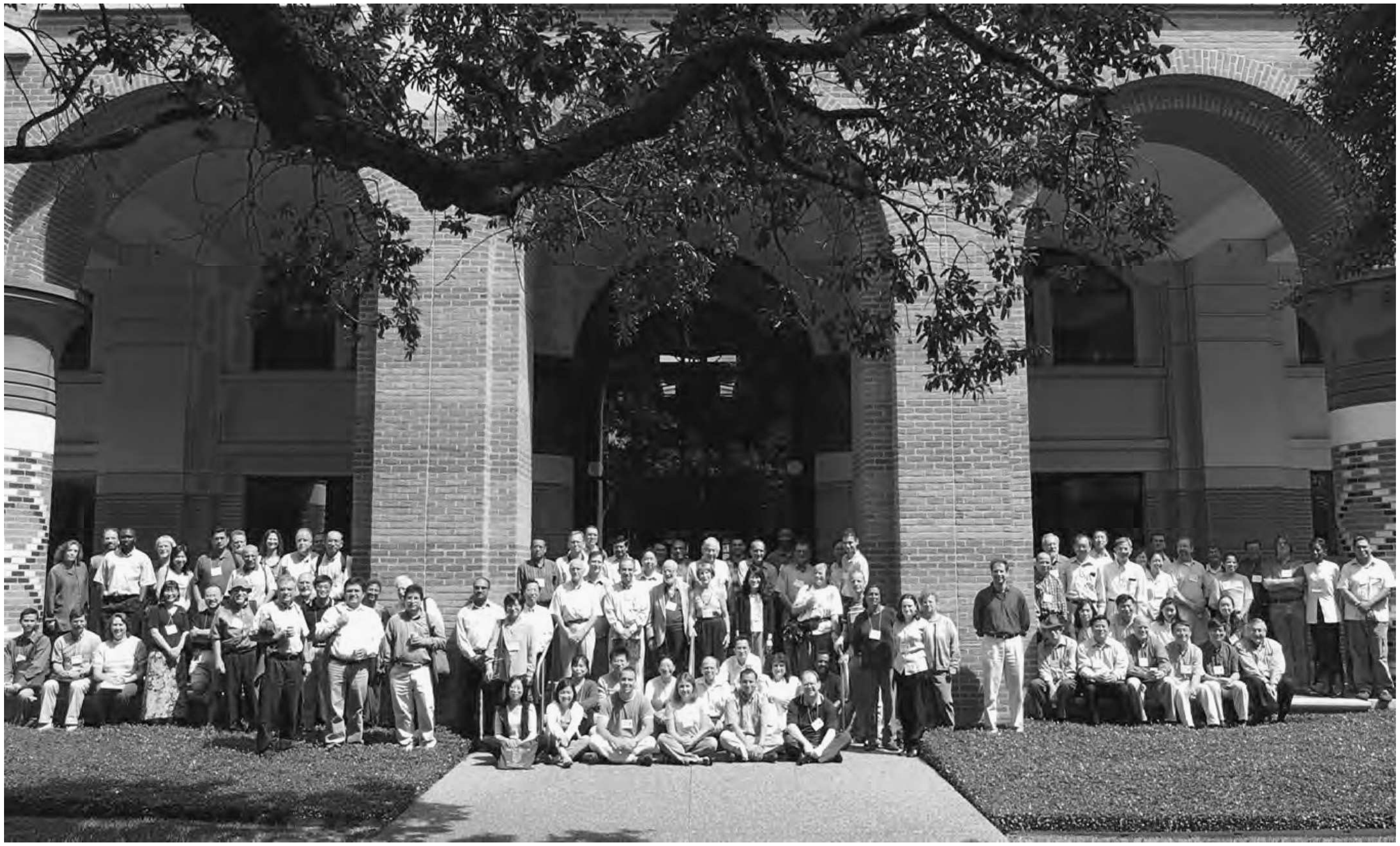

The Second Lehmann Symposium—Optimality
Rice University, May 19–22, 2004

## Partial List of Participants

Asheber Abebe
Auburn University
abebeas@auburn.edu

Hansjoerg Albrecher
Graz University
albrecher@tugraz.at

Demissie Alemayehu
Pfizer
alem@stat.columbia.edu

E. Neely Atkinson
University of Texas
MD Anderson Cancer Center
eatkinso@mdanderson.org

Scott Baggett
Rice University
baggett@rice.edu

Sarah Baraniuk
University of Texas Houston
School of Public Health
sbaraniuk@sph.uth.tmc.edu

Jose Luis Batun
CIMAT
batun@cimat.mx

Debasis Bhattacharya
Visva-Bharati, India
Debases_us@yahoo.com

Chad Bhatti
Rice University
bhatticr@rice.edu

Peter Bickel
University of California, Berkeley
bickel@stat.berkeley.edu

Sharad Borle
Rice University
sborle@rice.edu

Patrick Brockett
University of Texas, Austin
brockett@mail.utexas.edu

Barry Brown
University of Texas
MD Anderson Cancer Center
bwb@mdanderson.org

Ferry Butar Butar
Sam Houston State University
mth_fbb@shsu.edu

Raymond Carroll
Texas A&M University
carroll@stat.tamu.edu

Wenyaw Chan
University of Texas, Houston
Health Science Center
Wenyaw.Chan@uth.tmc.edu

Jamie Chatman
Rice University
jchatman@rice.edu

Cheng Cheng
St Jude Hospital
cheng.cheng@stjude.org

Hyemi Choi
Seoul National University
hyemichoi@yahoo.com

Blair Christian
Rice University
blairc@rice.edu

Daren B. H. Cline
Texas A&M University
dcline@stat.tamu.edu

Daniel Covarrubias
Rice University
dcorvarru@stat.rice.edu

David R. Cox
Nuffield College, Oxford
david.cox@nut.ox.ac.uk

Dennis Cox
Rice University
dcox@rice.edu

Kalatu Davies
Rice University
kdavies@rice.edu

Ginger Davis
Rice University
gmdavis@rice.edu

Richard Davis
Colorado State University
rdavis@stat.colostate.edu

Victor H. de la Peña
Columbia University
vp@stat.columbia.edu

Li Deng
Rice University
lident@rice.edu

Victor De Oliveira
University of Arkansas
vdo@uark.ed

Persi Diaconis
Stanford University

Chris Ding
Lawrence Berkeley Natl Lab
chqding@lbl.gov

Kjell Doksum
University of Wisconsin
doksum@stat.wisc.edu

Joan Dong
University of Texas
MD Anderson Cancer Center
qdong@mdanderson.org

Wesley Eddings
Kenyon College
eddingsw@kenyon.edu

Brad Efron
Stanford University
brad@statistics.stanford.edu

Hammou El Barmi
Baruch College
hammou_elbarmi@baruch.cuny.edu

Kathy Ensor
Rice University
kathy@rice.edu

Alan H. Feiveson
Johnson Space Center
alan.h.feiveson@nasa.gov

Hector Flores
Rice University
hflores@rice.edu

Garrett Fox
Rice University
gfox@stat.rice.edu

David A. Freedman
University of California, Berkeley
freedman@stat.berkeley.edu

Wenjiang Fu
Texas A&M University
wfu@stat.tamu.edu

Joseph Gastwirth
George Washington University
jlgast@gwu.edu

Susan Geller
Texas A&M University
geller@math.tamu.edu

Musie Ghebremichael
Rice University
musie@rice.edu

Irene Gijbels
Catholic University of Louvin
gijbels@stat.ucl.ac.be

Nancy Glenn
University of South Carolina
nglenn@stat.sc.edu

Carlos Gonzalez Universidad
Veracruzana
cglezand@tema.cum.mx

Shulamith Gross
NSF
sgross@nsf.gov

Xiangjun Gu
University of Texas
MD Anderson Cancer Center
xgu@mdanderson.org

Rudy Guerra
Rice University
rguerra@rice.edu

Shu Han
Rice University
shuhan@rice.edu

Robert Hardy
University of Texas
Health Science Center, Houston SPH
bhardy@sph.uth.tmc.edu

Jeffrey D. Hart
Texas A&M University
hart@stat.tamu.edu

Mike Hernandez
University of Texas
MD Anderson Cancer Center
Mike@sph.uth.tmc.edu

Richard Heydorn
NASA
richard.p.heydorn@nasa.gov

Tyson Holmes
Stanford University
tholmes@stanford.edu

Charlotte Hsieh
Rice University
hsiehc@rice.edu

Pushing Hsieh
University of California, Davis
fushing@wald.ucdavis.edu

Xuelin Huang
University of Texas
MD Anderson Cancer Center
xlhuang@mdanderson.org

Gabriel Huerta
University of New Mexico
ghuerta@stat.unm.edu

Sigfrido Iglesias
Gonzalez University of Toronto
sigfrido@fisher.utstat.toronto.c

Yuan Ji
University of Texas
yuanji@mdanderson.org

Sergio Juarez
Veracruz University
Mexico
sejuarez@uv.mx

Asha Seth Kapadia
University of Texas
Health Science Center, Houston SPH
School of Public Health
akapadia@sph.uth.tmc.edu

Vladislav Karguine
Cornerstone Research
slava@bu.edu

K. Krishnamoorthy
University of Louisiana
krishna@louisiana.edu

Mike Lecocke
Rice University
mlecocke@stat.rice.edu

Eun-Joo Lee
Texas Tech University
elee@math.ttu.edu

J. Jack Lee
University of Texas
MD Anderson Cancer Center
jjlee@mdanderson.org

Jaechoul Lee
Boise State University
jaechlee@math.biostate.edu

Jong Soo Lee
Rice University
jslee@rice.edu

Young Kyung Lee
Seoul National University
itsgirl@hanmail.net

Hannes Leeb
Yale University
hannes.leeb@yale.edu

Erich Lehmann
University of California,
Berkeley
shaffer@stat.berkeley.edu

Lei Lei
University of Texas
Health Science Center, SPH
llei@sph.uth.tmc.edu

Wei-Yin Loh
University of Wisconsin
loh@stat.wisc.edu

Yen-Peng Li
University of Texas, Houston
School of Public Health
yli@sph.uth.tmc.edu

Yisheng Li
University of Texas
MD Anderson Cancer Center
ysli@mdanderson.org

Simon Lunagomez
University of Texas
MD Anderson Cancer Center
slunago@mdanderson.org

Matthias Matheas
Rice University
matze@rice.edu

Deborah Mayo
Virginia Tech
mayod@vt.edu

Robert Mnatsakanov
West Virginia University
rmnatsak@stat.wvu.edu

Jeffrey Morris
University of Texas
MD Anderson Cancer Center
jeffmo@odin.mdacc.tmc.edu

Peter Mueller
University of Texas
MD Anderson Cancer Center
pmueller@mdanderson.org

Bin Nan
University of Michigan
bnan@umich.edu

Loki Natarajan
University of California,
San Diego
loki@math.ucsd.edu

E. Shannon Neeley
Rice University
sneeley@rice.edu

Josue Noyola-Martinez
Rice University
jcnm@rice.edu

Ingram Olkin
Stanford University
iolkin@stat.stanford.edu

Peter Olofsson
Rice University

Bernard Omolo
Texas Tech University
bomolo@math.ttu.edu

Richard C. Ott
Rice University
rott@rice.edu

Galen Papkov
Rice University
gpapkov@rice.edu

Byeong U Park
Seoul National University
bupark@stats.snu.ac.kr

Emanuel Parzen
Texas A&M University
eparzen@tamu.edu

Bo Peng
Rice University
bpeng@rice.edu

Kenneth Pietz
Department of Veteran Affairs
kpietz@bcm.tmc.edu

Kathy Prewitt
Arizona State University
kathryn.prewitt@asu.edu

Adrian Raftery
University of Washington
raftery@stat.washington.edu

Vinay Ribeiro
Rice University
vinay@rice.edu

Peter Richardson
Baylor College of Medicine
peterr@bcm.tmc.edu

Rolf Riedi
Rice University
riedi@rice.edu

Javier Rojo
Rice University
jrojo@rice.edu

Joseph Romano
Stanford University
romano@stat.stanford.edu

Gary L. Rosner
University of Texas
MD Anderson Cancer Center
glrosner@mdanderson.org

Mark Rothmann
US Food and Drug
Administration
rothmann@cder.fda.gov

Chris Rudnicki
Rice University
rudnicki@stat.rice.edu

Jose Aimer Sanqui
Appalachian St. University
sanquijat@appstate.edu

William R. Schucany
Southern Methodist University
schucany@smu.edu

Alena Scott
Rice University
oetting@stat.rice.edu

David W. Scott
Rice University
scottdw@rice.edu

Juliet Shaffer
University of California, Berkeley
shaffer@stat.berkeley.edu

Yu Shen
University of Texas
MD Anderson Cancer Center
yushen@mdanderson.org

Nozer Singpurwalla
The George Washington University
nozer@gwu.edu

Tumulesh Solanky
University of New Orleans
tsolanky@uno.edu

Julianne Souchek
Department of Veteran Affairs
jsoucheck@bcm.tmc.edu

Melissa Spann
Baylor University
melissa-spann@baylor.edu

Aris Spanos
Virginia Tech
aris@vt.edu

Hsiguang Sung
Rice University
hgsung@rice.edu

Gábor Székely
Bowling Green Sate
University
gabors@bgnet.bgsu.edu

Jef Teugels
Katholieke Univ. Leuven
jef.teugels@wis.kuleuven.be

James Thompson
Rice University
thomp@rice.edu

Jack Tubbs
Baylor University
jack_tubbs@baylor.edu

Jane-Ling Wang
University of California, Davis
wang@wald.ucdavis.edu

Naisyin Wang
Texas A&M University
nwang@stat.tamu.edu

Kyle Wathen
University of Texas
MD Anderson Cancer Center
& University of Texas GSBS
jkwathen@mdanderson.org

Peter Westfall
Texas Tech University
westfall@ba.ttu.edu

William Wojciechowski
Rice University
williamc@rice.edu

Jose-Miguel Yamal
Rice University &
University of Texas
MD Anderson Cancer Center
jmy@stat.rice.edu

Jun Yan
University of Iowa
jyan@stat.wiowa.edu

Guosheng Yin
University of Texas
MD Anderson Cancer Center
gyin@odin.mdacc.tmc.edu

Zhaoxia Yu
Rice University
yu@rice.edu

Issa Zakeri
Baylor College of Medicine
izakeri@bcm.tmc.edu

Qing Zhang
University of Texas
MD Anderson Cancer Center
qzhang@mdanderson.org

Hui Zhao
University of Texas
Health Science Center
School of Public Health
hzhao@sph.uth.tmc.edu

Yichum Zhao
Georgia State University
yzhao@math.stat.gsu.edu

*Acknowledgement of referees' services*

The efforts of the following referees are gratefully acknowledged

Jose Luis Batun
CIMAT, Mexico

Roger Berger
Arizona State University

Prabir Burman
University of California, Davis

Ray Carroll
Texas A&M University

Cheng Cheng
St. Jude's Children's Research Hospital

David R. Cox
Nuffield College, Oxford

Dorota M. Dabrowska
University of California, Los Angeles

Victor H. de la Pena
Columbia University

Kjell Doksum
University of Wisconsin, Madison

Armando Dominguez
CIMAT, Mexico

Sandrine Dudoit
University of California, Berkeley

Richard Dykstra
University of Iowa

Bradley Efron
Stanford University

Harnmou El Barmi
The City University of New York

Luis Enrique Figueroa
Purdue University

Joseph L. Gastwirth
George Washington University

Marc G. Genton
Texas A&M University

Musie Ghebremichael
Yale University

Graciela Gonzalez
Mexico

Hannes Leeb
Yale University

Erich L. Lehmann
University of California, Berkeley

Ker-Chau Li
University of California, Los Angeles

Wei-Yin Loh
University of Wisconsin, Madison

Hari Mukerjee
Wichita State University

Loki Natarajan
University of California, San Diego

Ingram Olkin
Stanford University

Liang Peng
Georgia Institute of Technology

Joseph P. Romano
Stanford University

Louise Ryan
Harvard University

Sanat Sarkar
Temple University

William R. Schucany
Southern Methodist University

David W. Scott
Rice University

Juliet P. Shaffer
University of California, Berkeley

Nozer D. Singpurwalla
George Washington University

David Sprott
CIMAT and University of Waterloo

Jef L. Teugels
Katholieke Universiteit Leuven

Martin J. Wainwright
University of California, Berkeley

Jane-Ling Wang
University of California, Davis

Peter Westfall
Texas Tech University

Grace Yang
University of Maryland, College Park

Yannis Yatracos (2)
National University of Singapore

Guosheng Yin
University of Texas MD Anderson Cancer Center

Hongyu Zhao
Yale University